\documentclass[10pt, dvips]{article}
\usepackage{epsfig}
\begin{document}

\newcommand{\p}{\partial}
\newcommand{\om}{\omega}
\newcommand{\e}{\epsilon}
\renewcommand{\a}{\alpha}
\renewcommand{\b}{\beta}


\title{Escape Probability, Mean Residence Time and \\ 
	Geophysical Fluid Particle Dynamics}  
	
\author{James R. Brannan, Jinqiao Duan\footnote{Author for correspondence --
Dr. Jinqiao Duan, E-mail:  duan@math.clemson.edu  $\;$ 
Fax: (864)656-2730 $\;$ Tel: (864)656-2730 } $\;$  and Vincent J. Ervin \\
\\
Clemson University \\
Department of Mathematical Sciences \\
Clemson, South Carolina 29634, USA.  }

\date{December 25, 1998 (Final Version) }

\maketitle

 
{\tiny

\begin{abstract}
 
Stochastic dynamical systems arise as models for
fluid particle motion  in geophysical flows with
random velocity fields. Escape probability (from a fluid domain) and mean
residence time (in a fluid domain)    quantify
fluid transport between flow
regimes of different characteristic motion.  

We consider a quasigeostrophic meandering jet model with random
perturbations. This jet is parameterized by the  
parameter $\beta =  \frac{2 \Omega}{r} \cos (\theta)$, where
$\Omega$ is the rotation rate of the earth, $r$ the earth's radius and $\theta$
the latitude. Note that  $\Omega$ and $r$ are fixed, so $\beta$ is a monotonic
decreasing function of the latitude.   The unperturbed jet
(for $0 < \beta <\frac23$)
consists of a basic flow   with attached eddies.
With random perturbations, there is
fluid exchange between regimes of different characteristic motion.
We quantify the exchange by escape probability and mean residence time. 

    For an eddy, the average  escape probability for fluid particles 
(initially inside the eddy) escape into the exterior retrograde region 
is smaller than escape into the jet core for $0 < \beta < 0.3333$, 
while for $0.3333  <\beta < \frac23$, the opposite holds. 
  
    For a unit jet core near the jet troughs, the average escape probability for
fluid particles (initially inside the 
jet core)  escape into the northern recirculating region is greater
than escape into the southern recirculating region for $0 < \beta < 0.115$, 
while for $0.385 < \beta < \frac23$, the opposite holds. 
Moreover, for $0.115 < \beta < 0.385$,
fluid particles are about equally likely to escape into either 
recirculating regions.  
 Furthermore, for a unit jet core near the jet crests,  
the situation is the opposite as for near the jet troughs.

The maximal mean residence time of fluid particles initially
in an eddy increases as $\beta$  increases  from $0$ to $0.432$
(or as latitude decreases accordingly),
then decreases as $\beta$ increases  from $0.432$ to $\frac23$
(or as latitude decreases accordingly).
However, the maximal mean residence time of fluid particles initially
in a unit jet core always increases as $\beta $ increases (or as
latitude decreases).


\end{abstract}

}

\newpage

\section{Stochastic dynamics: Escape probability and mean residence time}

Stochastic dynamical systems are used as models for various scientific and
engineering problems.
We consider the following class of stochastic dynamical systems
\begin{eqnarray}
\dot{x} & = & a_1(x,y) + b_1(x,y) \dot{w}_1, \label{eqn1} \\ 
\dot{y} & = & a_2(x,y) + b_2(x,y) \dot{w}_2, \label{eqn2}  
\end{eqnarray}
where  $w_1(t), w_2(t)$  
are two real independent Wiener processes (white noises)  and
$a_1, a_2, b_1, b_2$ are given deterministic functions.
More complicated stochastic systems also occur
in applications (\cite{Arnold}, \cite{Oksendal}, \cite{Freidlin}), 
\cite{Freidlin2}, \cite{Stone}, \cite{Kloeden}). 

For a planar bounded domain $D$, we can consider the exit problem
of random solution trajectories of (\ref{eqn1})-(\ref{eqn2}) 
from $D$. To this end, let $\p D$ denote the boundary of $D$
and let $\Gamma$ be a part of the boundary $\p D$.
The escape probability $p(x,y)$ is the probability   
that the trajectory of a particle starting at $(x,y)$ in $D$
first hits $\p D$ (or escapes from $D$) at some point in $\Gamma$, and 
$p(x,y)$ is known to satisfy
(\cite{Lin}, \cite{Schuss} and references therein)

\begin{eqnarray} 
  \frac12 b_1^2(x,y) p_{xx} + \frac12 b_2^2(x,y) p_{yy} + a_1(x,y) p_x
               +a_2(x,y) p_y  & = &  0,  \label{eqn3} \\
    p|_{\Gamma} & = & 1,         \label{eqn4}  \\
    p|_{\p D - \Gamma} & = & 0. \label{eqn5}
\end{eqnarray} 
Suppose that initial conditions (or initial particles) are 
uniformly distributed over $D$. 
The average escape probability $ P $ that a trajectory
will leave $D$ along the subboundary $\Gamma$, before leaving the
rest of the boundary, is given by (e.g.,  \cite{Lin}, \cite{Schuss})

\begin{equation}
	  P =  \frac{1}{|D|} \int\int_D p(x,y) dxdy,
	  \label{average}
\end{equation} 
where $|D|$ is the area of domain $D$. 
 
The residence time of a particle initially at $(x,y)$ inside $D$
is the time until the particle first hits $\p D$ (or escapes from $D$).
The mean residence time $u(x,y)$ is given by 
(e.g., \cite{Schuss}, \cite{Naeh}, \cite{Risken} and references therein)

\begin{eqnarray} 
    \frac12 b_1^2(x,y) u_{xx} + \frac12 b_2^2(x,y) u_{yy} + a_1(x,y) u_x
               +a_2(x,y) u_y  & = &  -1,   \label{eqn6}\\
 	         u|_{\p D}    & = &    0.  \label{eqn7}
\end{eqnarray}

\section{A quasigeostrophic jet model} 

The Lagrangian view of fluid motion is particularly important in geophysical
flows since only Lagrangian data can be obtained in many situations. 
It is essential to understand fluid particle trajectories
in many fluid problems. Escape probability (from a fluid domain) and mean
residence time (in a fluid domain)    quantify
fluid transport between flow
regimes of different characteristic motion. Deterministic quantities
like escape probability and  mean  residence time can be computed
by solving Fokker-Planck type partial differential equations. 

We now use the these ideas in the
investigation of meandering oceanic jets.
Meandering oceanic jets such as the Gulf Stream are strong
currents dividing different bodies of water.

Recently, del-castillo-Negrete and Morrison (\cite{Diego}), and
Pratt et al. (\cite{Pratt}, \cite{Miller}) have studied
models for oceanic jets. These models are
dynamically consistent to within a linear approximation,
i.e., the potential vorticity is approximately conserved.
Del-castillo-Negrete and Morrison's model consists
of the basic flow plus time-periodic linear neutral
modes.

In this paper, we consider an oceanic jet consisting of the basic flow
as in del-castillo-Negrete and Morrison (\cite{Diego}),
plus random-in-time noise. This model incorporates
small-scale oceanic motions such as the molecular
diffusion (\cite{Isichenko}), which is an important
factor in the Gulf Stream (\cite{Dutkiewicz}, \cite{Griffa2}).
The irregularity of RAFOS floats (\cite{Bower},
\cite{Song}, \cite{Lozier}) also suggests the inclusion of
random effects in Gulf Stream modeling.

This random jet may also be viewed as satisfying,
approximately in the spirit of del-castillo-Negrete and 
Morrison (\cite{Diego}), the randomly wind forced
quasigeostrophic model.
Several authors have considered the randomly wind forced
quasigeostrophic model in order to incorporate the impact of
uncertain geophysical forces (\cite{Samelson}, 
\cite{Griffa}, \cite{Holloway}, \cite{Muller},
\cite{DelSole-Farrell}, \cite{Duan}).  
They studied statistical issues such as
estimating correlation coefficients for the {\em linearized}
quasigeostrophic equation with random forcing. There is also recent
work which investigates the   impact of the uncertainty of
the ocean bottom topography on
quasigeostrophic dynamics (\cite{Klyatskin_Gurarie}).

The randomly forced quasigeostrophic equation takes the form
(\cite{Muller})
\begin{equation}
   \Delta \psi_t + J(\psi, \Delta \psi ) + \beta \psi_x
 		  =  \frac{dW}{dt} \; ,
   \label{qg}
\end{equation}
where $W(x,y,t)$ is a space-time Wiener process (white noise).
The stream function would have  a random   or noise component
(\cite{Griffa2}, \cite{Dutkiewicz}).
Note that $\beta$ is the meridional derivative of the Coriolis
parameter \cite{Pedlosky} \cite{Benoit}, i.e., 
$\beta =  \frac{2 \Omega}{r} \cos (\theta)$, where
$\Omega$ is the rotation rate of the earth, $r$ the earth's radius and $\theta$
  the latitude. Since  $\Omega$ and $r$ are fixed, $\beta$ is a monotonic
  decreasing function of the latitude.

The deterministic meandering jet
derived in  del-castillo-Negrete and Morrison (\cite{Diego}) is 
$$
\Psi(x,y) = -\tanh (y) + a \; sech^2 (y) \cos(kx) + cy,
$$
where 
$$
a =0.01, c=\frac13 (1+\sqrt{1-\frac32 \beta}),
k=\sqrt{2(1+\sqrt{1-\frac32 \beta})},  0 \leq \beta \leq \frac23.
$$
This $\Psi(x,y)$ is an approximate solution
of the usual quasigeostrophic model
\begin{equation}
   \Delta \psi_t + J(\psi, \Delta \psi ) + \beta \psi_x
 		  = 0 \; .
   \label{usual-qg}
\end{equation}

With random  wind forcing  or molecular diffusive forcing
in the stochastic quasigeostrophic model
(\ref{qg}),  the stream function would have a 
random or noise component. We approximate this  noise component
by adding a noise term to the above deterministic stream function
$\Psi(x,y)$, that is, in the rest of this paper, we consider the
following random stream function as a model for a
quasigeostrophic meandering jet,
$$
\tilde{\Psi} (x,y) = -\tanh (y) + a \; sech^2 (y) \cos(kx) + cy + \mbox{noise}.
$$

The equations of motion for fluid particles in this   jet then
have noise terms. We further approximate them 
as white noises (or Wiener processes)  
\begin{eqnarray}
dx & = & -\Psi_y   dt + \sqrt{\e} dw_1, \\ 
dy & = & \;\; \Psi_x  dt  + \sqrt{\e} dw_2,  
\end{eqnarray}
or more specifically,
\begin{eqnarray}
dx & = & [sech^2(y) + 2a \;sech^2(y) \tanh (y) \cos(kx) -c] dt 
		+ \sqrt{\e} dw_1,  \label{finaleqn1}  \\ 
dy & = &  -ak\; sech^2(y) \sin(kx)  dt 
		+  \sqrt{\e} dw_2,  \label{finaleqn2}
\end{eqnarray}
where $0 < \e < 1$, and $w_1(t), w_2(t)$  
are two real independent Wiener processes (in time only).
The calculations below are for $\e = 0.001 $.
Note that $\beta$ is now the only parameter \ref{finaleqn1}, \ref{finaleqn2}, 
as $a$ is given and $c,$ and $ k$ depend only on $\beta \in [0, \frac23]$.

\begin{figure}[tbp] 
\hbox to \hsize{\hfil \psfig{figure=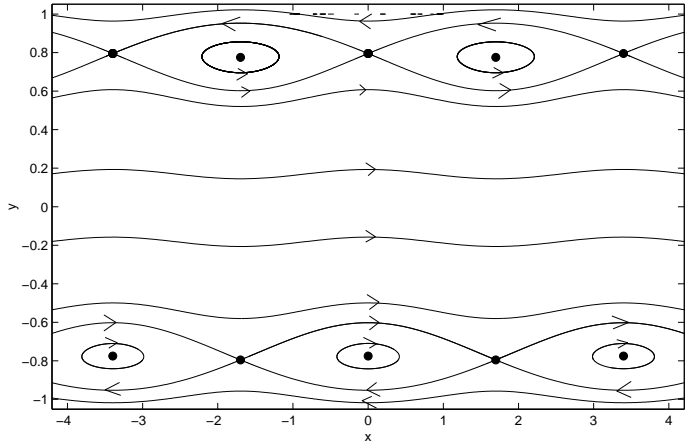,height=3in}\hfil}
\caption{Unperturbed jet: $\e =0$ and $\beta=1/3$. }
\label{unperturb}
\end{figure}

When $\e = 0$, the deterministic jet
consists of the jet core and two rows of
recirculating eddies, which are called the northern and southern
recirculating regions. Outside the recirculating regions are the
exterior retrograde regions; see Figure \ref{unperturb}.

\section{Escape probability} 

\begin{figure}[tp]
\hbox to \hsize{\hfil \psfig{figure=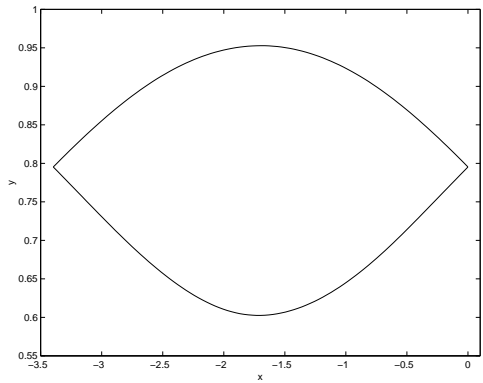,height=2in}\hfil}
\caption{An eddy: $\beta=1/3$ }
\label{eddy}
\end{figure}

\begin{figure}[htbp] 
\hbox to \hsize{\hfil \psfig{figure=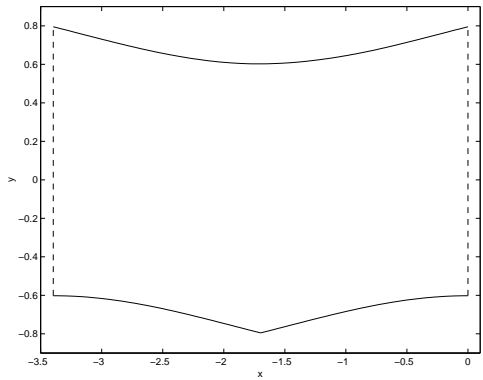,height=2in}\hfil}
\caption{A  unit jet core near a trough: $\beta=1/3$}
\label{jet}
\end{figure}

We take $D$ to be either an eddy or
a piece of jet core (see Figures \ref{eddy}, \ref{jet}).
This piece of jet core has the same horizontal length scale as an eddy, and
it is one period of the deterministic jet core (note that the deterministic
velocity field is periodic in $x$). Thus we call this piece of jet core 
a unit jet core.

From (\ref{eqn3}), (\ref{eqn4}), (\ref{eqn5}), 
the escape probability $p(x,y)$ that a 
fluid particle, initially at $(x,y)$,
crosses the subboundary $\Gamma$ of  the domain $D$ satisfies
 
\begin{eqnarray} 
  \e \Delta p + [sech^2(y) + 2a \;sech^2(y) \tanh (y) \cos(kx) -c] p_x
               -ak\; sech^2(y) \sin(kx) p_y  & = &  0, \\
    p|_{\Gamma} & = & 1,          \\
    p|_{\p D - \Gamma} & = & 0.
\end{eqnarray}
We take $\Gamma$ to be either top or bottom boundary of
an eddy or a unit jet core (see Figures \ref{eddy}, \ref{jet}). 
We numerically solve this elliptic system for 
various values of $\beta$ between $0$ and $\frac23$.
In the unit jet core case, we take periodic boundary
condition in horizontal (meridional) $x$ direction,
with period $\frac{2\pi}{k}$.

\begin{figure}[tbp] 
\hbox to \hsize{\hfil \psfig{figure=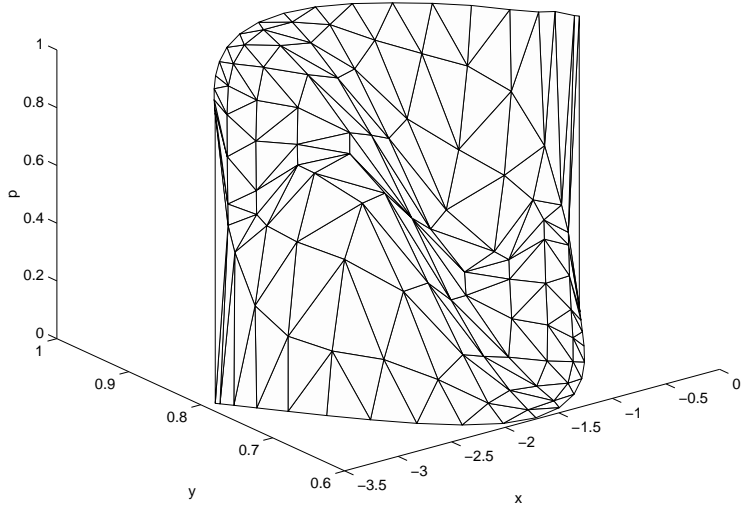,height=3in}\hfil}
\caption{Escape probability of fluid particles (initially in an eddy
in Figure \ref{eddy}) exiting into the exterior retrograde   
region: $\beta=1/3$  }
\label{up}
\end{figure}

\begin{figure}[tbp] 
\hbox to \hsize{\hfil \psfig{figure=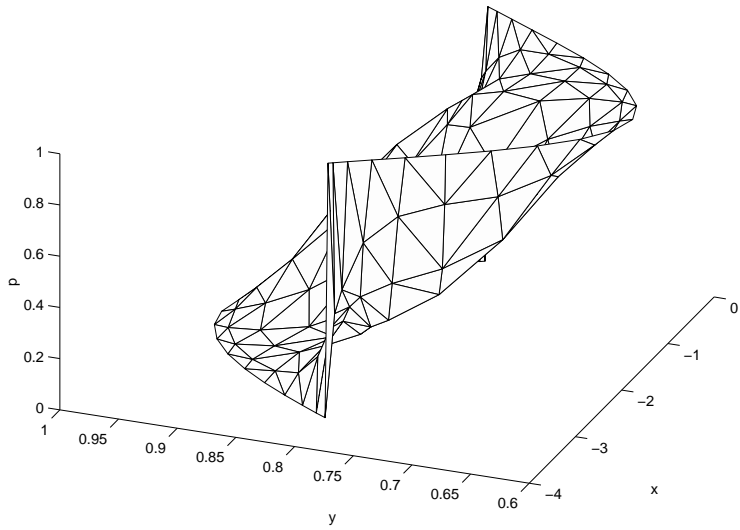,height=3in}\hfil}
\caption{Escape probability of fluid particles (initially in an eddy
in Figure \ref{eddy}) exiting into the jet core: $\beta=1/3$  }
\label{low}
\end{figure}

\begin{figure}[tbp] 
\hbox to \hsize{\hfil \psfig{figure=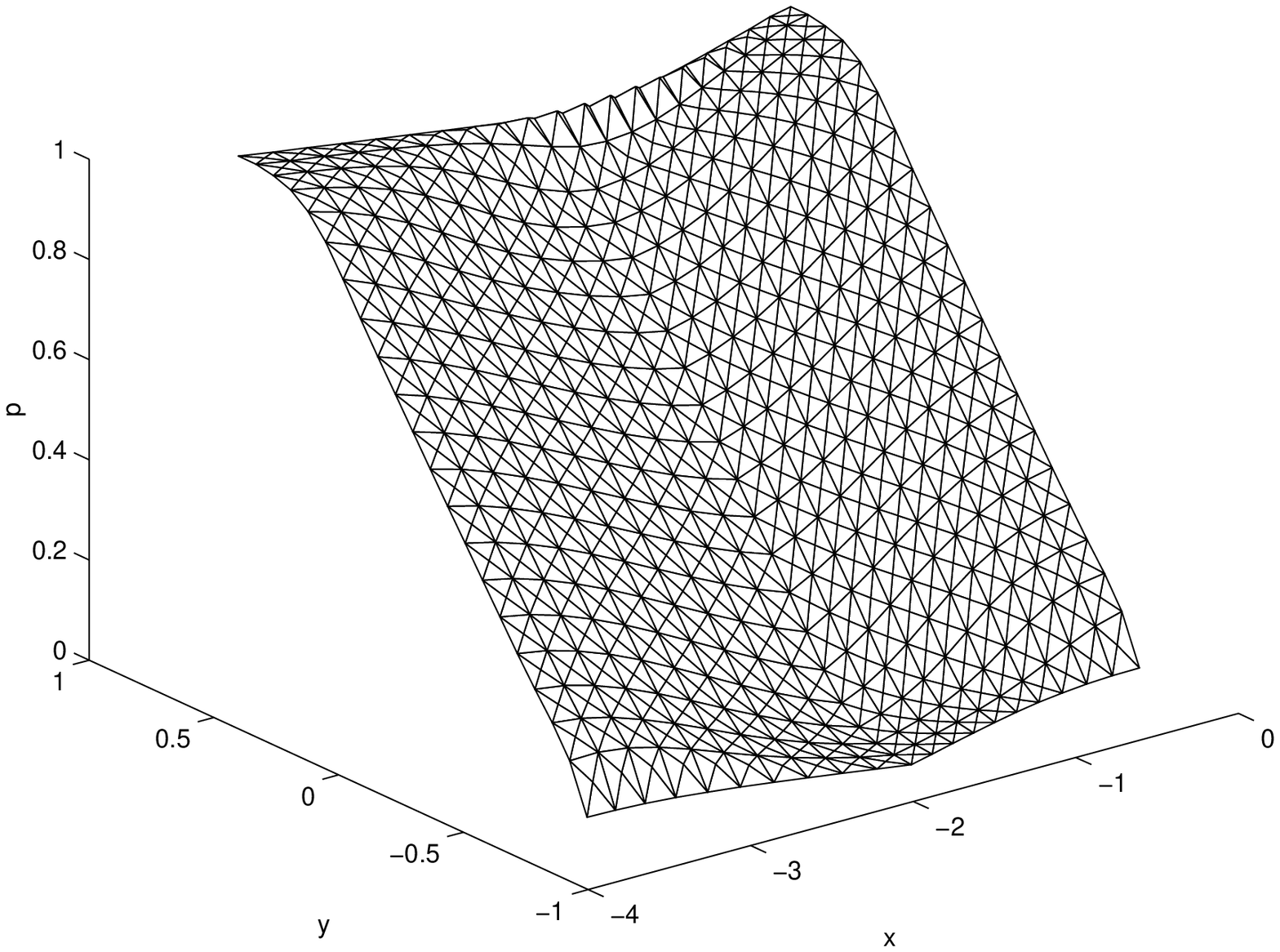,height=3in}\hfil}
\caption{Escape probability of fluid particles (initially in the unit jet core
in Figure \ref{jet}) exiting into the northern recirculating 
region: $\beta=1/3$ }
\label{top}
\end{figure}

\begin{figure}[tbp] 
\hbox to \hsize{\hfil \psfig{figure=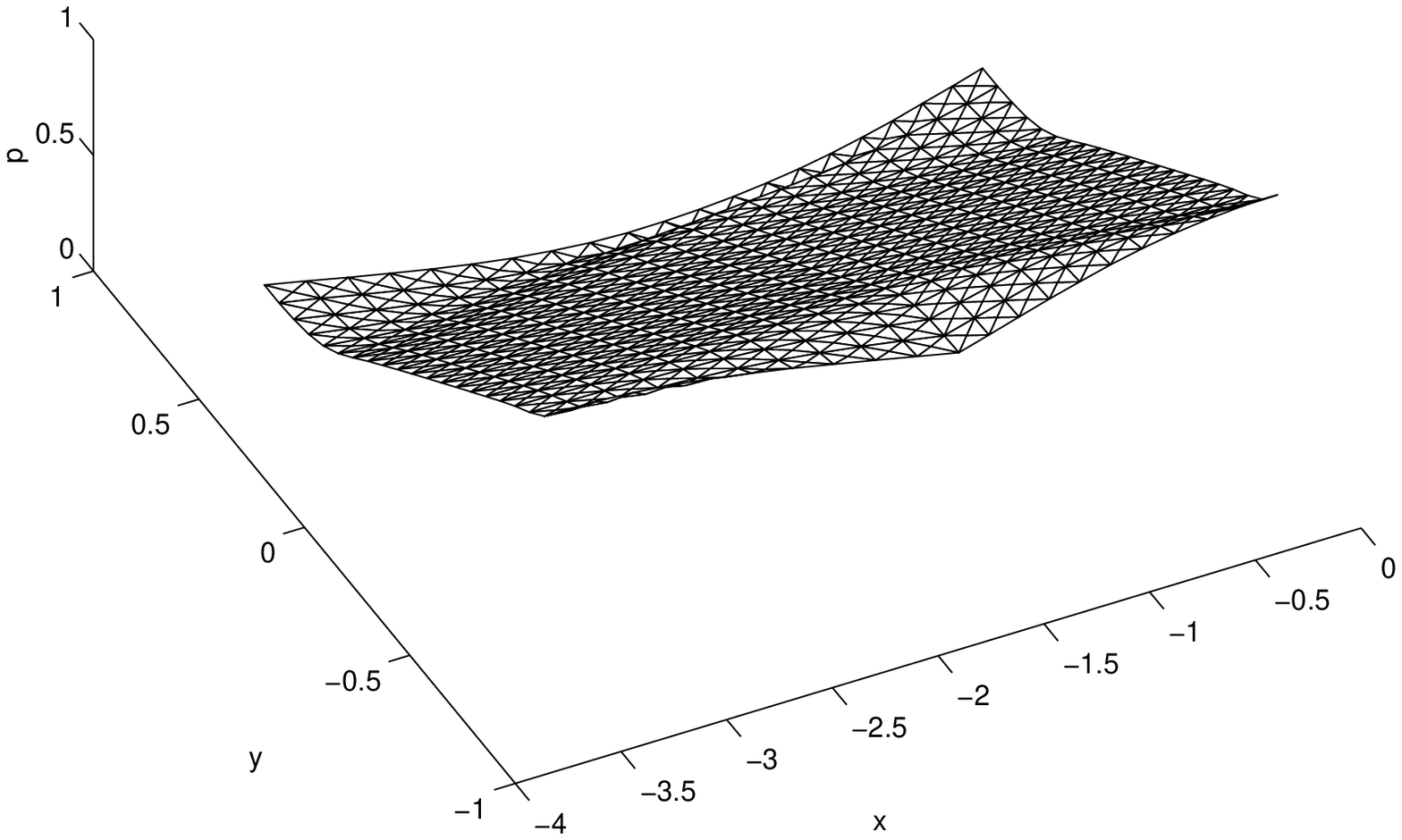,height=3in}\hfil}
\caption{Escape probability of fluid particles (initially in the unit jet core
in Figure \ref{jet}) exiting into the southern recirculating 
region: $\beta=1/3$  }
\label{bottom}
\end{figure}

A piecewise linear, finite element approximation scheme was used
for the numerical solutions of the escape probability $p(x,y)$,
and the mean residence time $u(x,y)$, described by the elliptic
equations (\ref{eqn3}), and (\ref{eqn6}), respectively. Using a 
collection of points lying on the boundary, piecewise cubic
splines were constructed to define the boundary of the eddy, and
the top and bottom boundaries of the jet core. Computational
(triangular) grids for the eddy and the jet core were then
obtained by deforming regular grids constructed for an ellipse
and a rectangular region, respectively. The computed escape probability
crossing the upper or lower boundary of an eddy or a unit jet core,
for the case of $\beta=\frac13$,
are shown in Figures \ref{up}, \ref{low}, \ref{top} and \ref{bottom}.

\begin{figure}[tbp] 
\hbox to \hsize{\hfil \psfig{figure=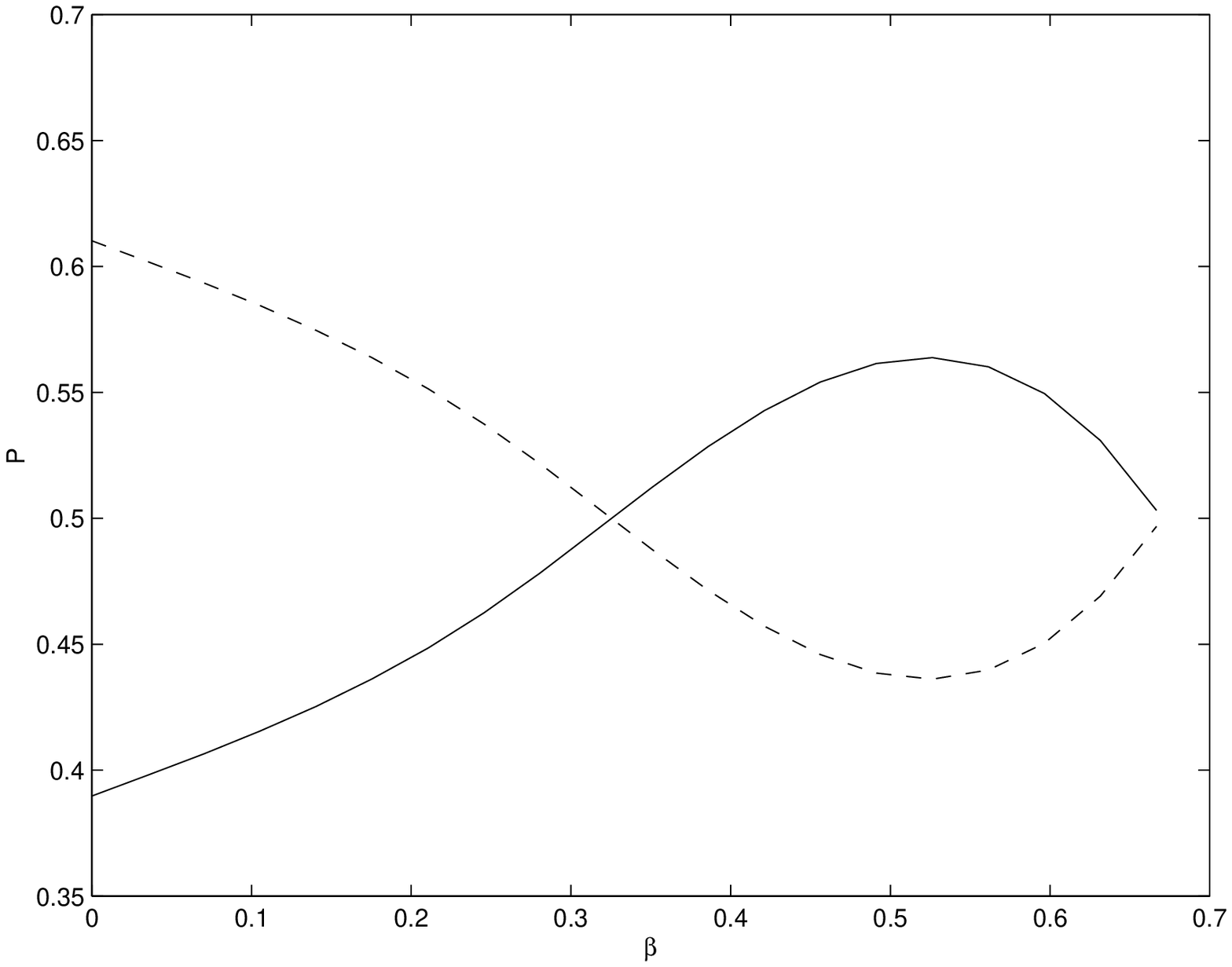,height=3in}\hfil}
\caption{Average escape probability for a fluid particle in an eddy:
  --- for exiting into the 
exterior retrograde region; $\_ \_ \_$ for exiting into the jet core region. }
\label{eddyaverage}
\end{figure}

\begin{figure}[tbp] 
\hbox to \hsize{\hfil \psfig{figure=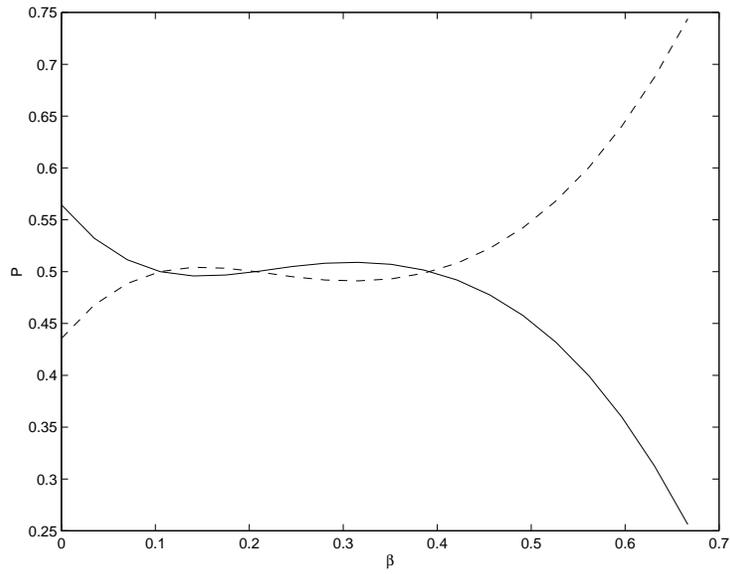,height=3in}\hfil}
\caption{Average escape probability for a fluid particle in the unit jet
core near a trough: --- for exiting into the northern 
recirculating region;  $\_ \_ \_$ for exiting into the southern 
recirculating region. }
\label{jetaverage}
\end{figure}

Suppose that the fluid particles are initially uniformly distributed
in $D$ (an eddy or a unit jet core). We can also compute the average escape probability 
$ P $ that a particle
will leave $D$ along the upper or lower subboundary $\Gamma$, using the
formula (\ref{average}); see Figures \ref{eddyaverage} and \ref{jetaverage}.
  
    For an eddy, the average  escape probability for fluid particles 
(initially inside the eddy) escape into the exterior retrograde region 
is smaller than escape into the jet core for $0 < \beta < 0.3333$, 
while for $0.3333  <\beta < \frac23$, the opposite holds 
(Figure \ref{eddyaverage}). Thus $\beta = 0.3333$ is
  a bifurcation point. Also, the average  escape probability
 for fluid particles  escape into the exterior retrograde region
 increases as $\beta$ increases from $0$ to $0.54$
 (or  as latitude decreases accordingly), and then 
 decreases as $\beta$ increases from $0.54$ to $\frac23$
 (or  as latitude decreases accordingly). Thus $\beta = 0.54$ is
  another bifurcation point. The opposite holds for
  the average  escape probability
 for fluid particles  escape into the jet core.
 
    For a unit jet core near the jet troughs, the average escape probability for
fluid particles (initially inside the 
jet core)  escape into the northern recirculating region is greater
than escape into the southern recirculating region for $0 < \beta < 0.115$, 
while for $0.385 < \beta < \frac23$, the opposite holds. 
Moreover, for $0.115 < \beta < 0.385$,
fluid particles are about equally likely to escape into either 
recirculating regions (Figure \ref{jetaverage}).  
 
Furthermore, for a unit jet core near the jet crests,  
the situation is the opposite as for near the jet troughs.

 \section{Mean residence time}

The mean residence time $u(x,y)$ of a fluid particle, initially at $(x,y)$ 
in either an eddy or
a piece of jet core (see Figures \ref{eddy}, \ref{jet}),  satisfies  
\begin{eqnarray} 
 \e \Delta u + [sech^2(y) + 2a \;sech^2(y) \tanh (y) \cos(kx) -c] u_x
    -ak\; sech^2(y) \sin(kx) u_y  & =  &  -1, 	\\
 	  u|_{\p D}    & = &    0
\end{eqnarray}

\begin{figure}[bp] 
\hbox to \hsize{\hfil \psfig{figure=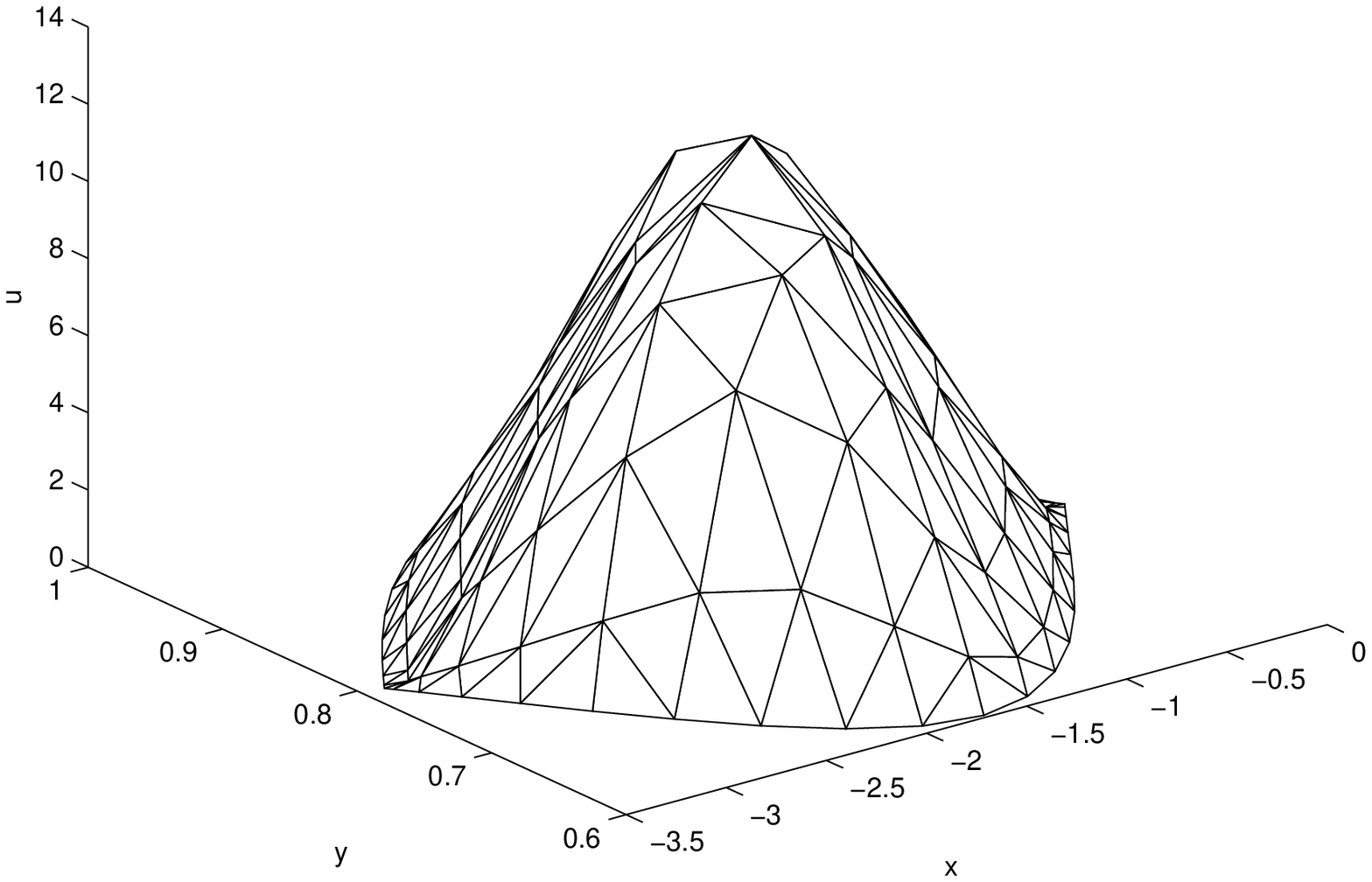,height=3in}\hfil}
\caption{Mean residence time in an eddy: $\beta = \frac13$ }
\label{eyemean}
\end{figure}

\begin{figure}[tbp] 
\hbox to \hsize{\hfil \psfig{figure=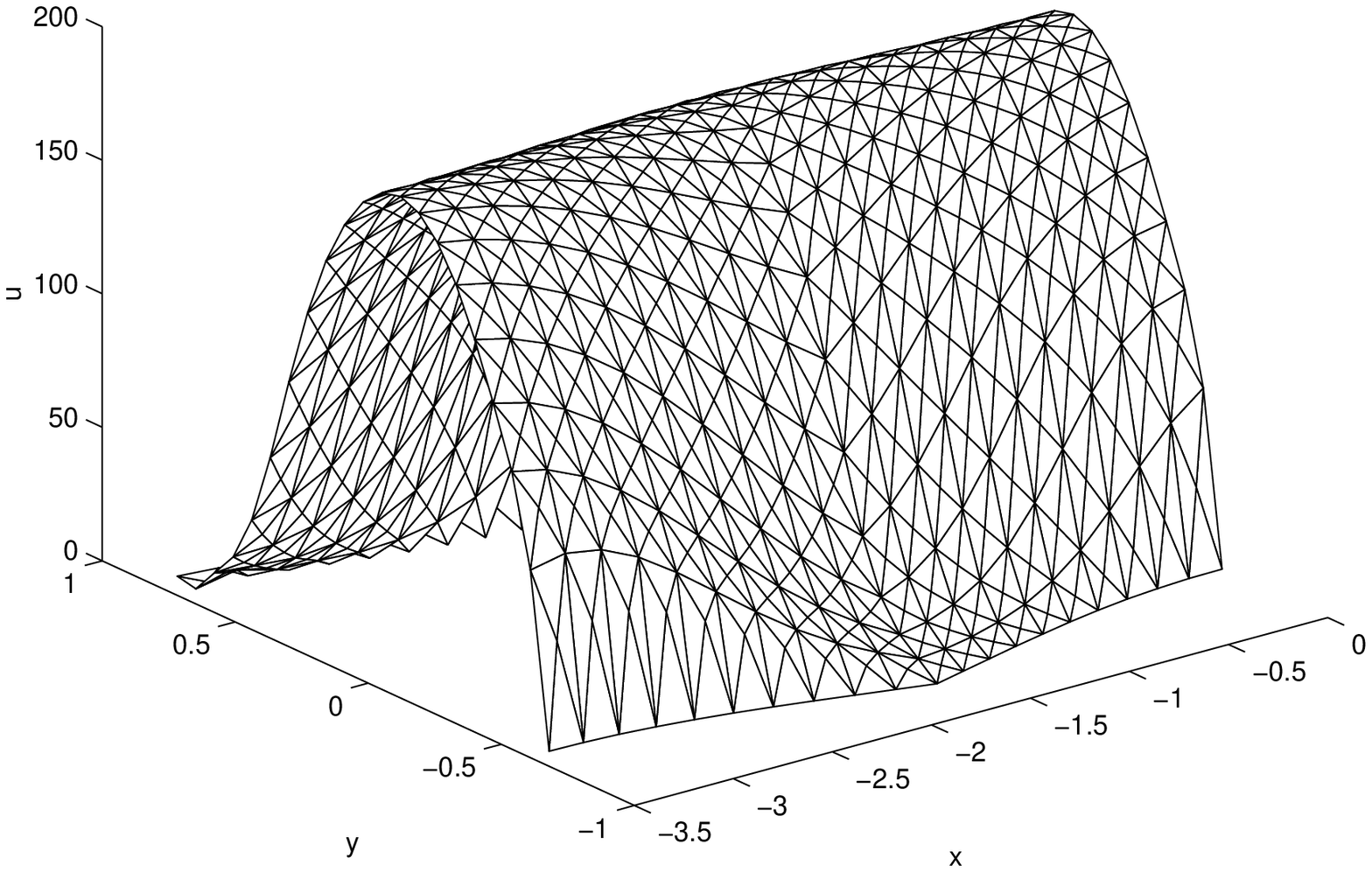,height=3in}\hfil}
\caption{Mean residence time in a  unit jet core: $\beta = \frac13$}
\label{jetmean}
\end{figure}

The mean residence times of fluid particles in an eddy or
a unit jet core are shown, for the case of $\beta =\frac13$,
 in Figures  \ref{eyemean}, \ref{jetmean}.

\begin{figure}[tbp] 
\hbox to \hsize{\hfil \psfig{figure=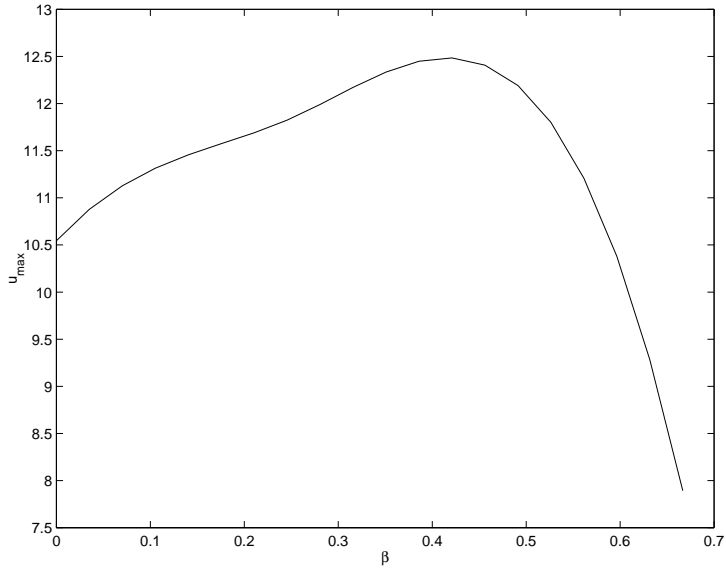,height=3in}\hfil}
\caption{Maximal value of mean residence time, as a function of $\beta$,  
in an eddy. }
\label{eyemax}
\end{figure}

\begin{figure}[tbp] 
\hbox to \hsize{\hfil \psfig{figure=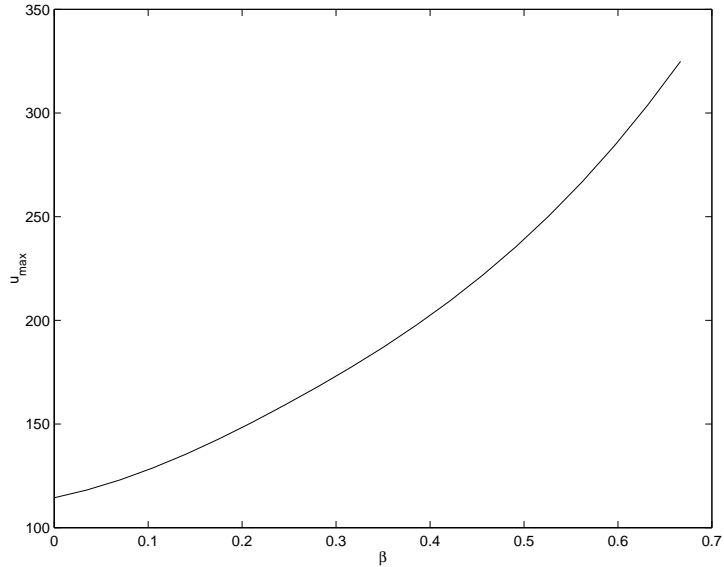,height=3in}\hfil}
\caption{Maximal value of mean residence time, as a function of $\beta$, 
in a unit jet core. }
\label{jetmax}
\end{figure}

The maximal mean residence time of fluid particles initially
in an eddy increases as $\beta$  increases  from $0$ to $0.432$
(or as latitude decreases accordingly),
then decreases  as $\beta$  increases from $0.432$ to $\frac23$
(or as latitude decreases accordingly);
see Figure \ref{eyemax}.
However, the maximal mean residence time of fluid particles initially
in a unit jet core always increases as $\beta $ increases
(or as latitude decreases accordingly); 
see Figure \ref{jetmax}.

\section{Discussions}

The present work on fluid particle motion in random flows
takes into account of fluid particle {\em diffusive} as well as advective motion.
There has been recent work
on fluid particle {\em advective} motion  (molecular diffusion ignored)
in time-periodic, 
quasi-periodic and aperiodic flows 
(periodic $\rightarrow$ quasi-periodic  $\rightarrow$ aperiodic $\rightarrow$
random); see, for example, \cite{Rom-Kedar}, \cite{Wiggins},
\cite{Samelson2}, \cite{Beigie}, \cite{Duan-Wiggins}, \cite{Wiggins2}
and \cite{Miller}.

Our work on random particle motion does {\em not} 
 require that the random part is small. However it does
require that the flow can be decomposed into steady or unsteady 
deterministic (drift)
and random (diffusion) parts; Otherwise the Fokker-Planck formalism
does not hold.

\bigskip

\noindent  {\bf Acknowledgement.} 
J. Duan would like to acknowledge the hospitality of 
the Center for Nonlinear Studies, Los Alamos National
Laboratory, the
Institute for Mathematics and its Applications (IMA), Minnesota,
and the  Woods Hole Oceanographic Institution, Massachusetts.
He is also grateful for the discussions on this work with
Diego del-castillo-Negrete (Scripps Institution of Oceanography),
Greg Holloway (Institute of Ocean Sciences, Canada),
Julian  Hunt (Cambridge University),  
Peter Kiessler (Clemson University),
Pat Miller (Stevens Institute of Technology),   
Larry Pratt (Woods Hole Oceanographic Institution),
and  Roger Samelson (Oregan State University).

    This work was supported by the National Science 
Foundation Grant DMS-9704345.

\end{document}